\documentclass[11pt, reqno]{amsart}
\usepackage[dvipsnames]{xcolor}
\usepackage{mathtools}
\usepackage{etoolbox}

\definecolor{mypink1}{rgb}{0.858, 0.188, 0.478}
\definecolor{mypink2}{RGB}{219, 48, 122}
\definecolor{mypink3}{cmyk}{0, 0.7808, 0.4429, 0.1412}
\definecolor{mygray}{gray}{0.6}

\definecolor{gris75}{gray}{0.25}
\definecolor{violet}{rgb}{0.5,0,0.5}

\definecolor{BrickRed}{rgb}{0.58, 0.0, 0.83}

\definecolor{armygreen}{rgb}{0.29, 0.33, 0.13}

\definecolor{brass}{rgb}{0.71, 0.65, 0.26}
\definecolor{antiquefuchsia}{rgb}{0.57, 0.36, 0.51}
\definecolor{amethyst}{rgb}{0.6, 0.4, 0.8}

\definecolor{mauvetaupe}{rgb}{0.57, 0.37, 0.43}

\usepackage{amsmath}
\usepackage{amssymb}
\usepackage{amsfonts}
\usepackage{mathrsfs}
\usepackage{amsbsy}
\usepackage{relsize}
\usepackage[colorlinks, linkcolor=LimeGreen, citecolor=blue, urlcolor=blue, hypertexnames=false]{hyperref}

\usepackage{tikz}

\numberwithin{equation}{section}

\numberwithin{bbb}{section}

\newtheorem{th1}{{\bf Theorem}}[section]
\newtheorem{thm}[th1]{{\bf Theorem}}
\newtheorem{lem}[th1]{{\bf Lemma}}

\newcommand{\R}{\mathbb{R}}

\newcommand{\C}{\mathbb{C}}

\usepackage{comment}
\PassOptionsToPackage{reqno}{amsmath}

\author[B. Ayed. Sabria \& T. Saanouni]{Sabria Ben Ayed \& Tarek Saanouni$^*$}
\thanks{* Corresponding author.}

\address[T. Saanouni]{Department of Mathematics, College of Science, Qassim University, Buraydah, Kingdom of Saudi Arabia.}
\email{\sl\color{blue}{t.saanouni@qu.edu.sa}}

\address[B. Ayed. Sabria]{Department of Mathematics, College of Science, Qassim University, Buraydah, Kingdom of Saudi Arabia.}
\email{\sl\color{blue}{S.ayad@qu.edu.sa}}


\title[INLS]{Sharp threshold of scattering versus non-scattering for some mass-sub-critical inhomogeneous NLS}

\date{\today}
\begin{document}
\begin{abstract}
This work investigates the long time asymptotic behavior of some inhomogeneous non-linear Schr\"odinger type equations. We give sharp a threshold of scattering versus non-scattering of mass solutions, depending on the source term. This work complements the literature to the inhomogeneous regime.
\end{abstract}
\hrule 
\maketitle
\hrule 

\thispagestyle{empty}

\section{Introduction}
Our primary object of interest in this paper is to investigates the initial valued problem for the inhomogeneous non-linear Schr\"odinger equation
\begin{equation}
\left\{
\begin{array}{ll}
    {\textnormal i}\partial _tu +\Delta u=|u|^{q-1}|x|^{-{b}}u;\\
u(0,\cdot)=u_{0},
\label{S}\tag{INLS}
\end{array}
\right.
\end{equation}

and its close cousin the Cauchy problem for the inhomogeneous generalized Hartree equation 

\begin{equation}
\left\{
\begin{array}{ll}
{\textnormal i}\partial_tw +\Delta w=|w|^{q-2}(J_\alpha *|\cdot|^{-b}|w|^q)|x|^{-b}w ;\\
w(0,\cdot)=w_0.
\label{H}\tag{INLH}
\end{array}
\right.
\end{equation}

Here, the wave functions are $u:=u(t,x)\in\C,w:=w(t,x)\in\C$, where the time-space variable is $(t,x)\in\mathbb R\times\mathbb R^N$. We call inhomogeneous term, the singular quantity $|\cdot|^{-{b}}$, where ${b}>0$. The Riesz-potential is defined on $\R^N$ by 
$$J_\alpha:=C_{N,\alpha}\,|\cdot|^{\alpha-N},\quad C_{N,\alpha}:=\frac{\Gamma(\frac{N-\alpha}2)}{\Gamma(\frac\alpha2)\pi^\frac{N}22^\alpha},\quad  0<\alpha<N.$$
Here and hereafter, we make the assumption, needed for the local existence of solutions \cite{tt},

\begin{align}\label{cnd}
\min\{b,\alpha,N-\alpha,N-b,1-2b+\alpha\}>0.
\end{align}

 The inhomogeneous nonlinear equation of Schr\"odinger type \eqref{S} models propagation of the beam in nonlinear optics and plasma physics. Indeed, stable high-power propagation in a plasma could be achieved by sending a precursor laser beam to create a channel with reduced electron density. This, in turn, lowers the nonlinearity within the channel \cite{clvt,tsg}. The inhomogeneous nonlinear equation of Hartree type \eqref{H} describes various physical phenomena. Specific instances of this model emerge in the mean-field limit of large systems of non-relativistic atoms and molecules, as well as in the propagation of electromagnetic waves in plasmas, among other applications \cite{dr,lbac,jfel}.\\

The existence of energy subcritical solutions to \eqref{S} was first established in \cite{gs}. This result was later revisited in \cite{cmg}, where solutions in Strichartz spaces were studied under additional conditions for \( N = 2, 3 \). The distinction between global existence and scattering versus finite-time blow-up below the ground state threshold was explored in \cite{fr1, fr2, fr3} through the concentration-compactness argument introduced by Kenig and Merle \cite{km}. This work was further extended in \cite{lca,cc} using the Dodson-Murphy method \cite{dm}, and the assumption of spherical symmetry was relaxed in \cite{cfgm}. Additional discussions on more general inhomogeneous terms can be found in \cite{at,dms}. See also \cite{mgjm,ylis,kls,chl,yckl,ak,agt} for the critical regime. The local/global well-posedness of \eqref{H} was first established in \cite{tt} using an adapted Gagliardo-Nirenberg type identity. Then, in \cite{sx}, the scattering theory was studied using the approach from \cite{dm} under the spherically symmetric assumption. More recently, \cite{bx} proved the scattering theory for the non-radial case, inspired by \cite{jm}. See also \cite{kls2,mgcx} for the critical regime.\\

Most of the works dealing with the long time asymptotic behavior of solutions to the Schr\"odinger problems treat the mass-super-critical regime. Conversely, the mass-sub-critical regime is less investigated despite that every datum with finite mass gives a global solution. It seems that the first works in this direction are \cite{gv1,gv2}, where the scattering was proved in the conformal space for $\alpha(N)<q-1<\frac{4}{N}$, where the Strauss exponent is $\alpha(N):=\frac{2-N+\sqrt{16N+(N-2)^2}}{4N}$. Then, the scattering in $L^2$ was established for the range $\frac{2}{N}<q-1<\frac4N$ in \cite{ytky,kn}. Note that the datum was assumed to be in some weighted space because of some breakdown of regularity of scattering \cite{gel,gel2}. For $q-1\leq\frac2N$, the scattering in $L^2$ can only occur for the zero solution, see \cite{barab,rtg,ws}. As for the mass-critical case $q-1=\frac4N$, the global well-posedness and
scattering in $L^2$ are known due to the list of works \cite{d1,d2,d3}. \\

It is the aim of this note to investigate the long time asymptotic behavior for the inhomogeneous Schr\"odinger type equations \eqref{S} and \eqref{H}. Indeed, we give a threshold of scattering versus non-scattering in $L^2$ of global mass-sub-critical defocusing solutions. The threshold depends on the source term exponent, namely there is no scattering in $L^2$ for low source term exponents. The novelty here is to treat the inhomogeneous case $b\neq0$. Indeed, the present work complements \cite{barab,ytky,sxa} to the inhomogeneous regime. To the authors knowledge, this work seems to be the first one dealing with the sharp threshold of $L^2$ scattering versus non-scattering of the inhomogeneous non-linear Schr\"odinger problems \eqref{S} and \eqref{H}. So, this note aims to fill in this gap in the literature of inhomogeneous NLS, since we are in a curious situation where we know that all solutions to \eqref{S} and \eqref{H} with a finite mass initial data are global, but
we have no understanding of the asymptotic behavior of a vast majority of such solutions.\\

 The rest of the paper is organized as follows. Section \ref{sec2} gives the main result and some useful estimates. Section \ref{sec3} is devoted to prove the non-scattering in $L^2$. Section \ref{sec4} establishes the $L^2$ scattering of global solutions. \\

For short, the Lebesgue spaces endowed with usual norms are denoted by

\begin{align*}
L^r&:=L^r({\mathbb R}^N),\quad \|\,\cdot\,\|_r:=\|\,\cdot\,\|_{L^r},\quad\|\,\cdot\,\|:=\|\,\cdot\,\|_2.
\end{align*}

We denote by $\Omega_R$, the centered ball of $\R^N$ with radius $R>0$ and $\Omega:=\Omega_1$. Finally, we let $(T^-,T^+)$ the maximal existence interval of an eventual solution to \eqref{S} or \eqref{H}.

\section{Background and Main results}\label{sec2}
This section contains the main contribution of this note and some useful standard estimates.
\subsection{Preliminary}
Let us denote the free Schr\"odinger kernel

\begin{align}
e^{{\textnormal i}t\Delta}u:= \mathcal F^{-1}\Big(e^{-{\textnormal i}t|\cdot|^2}\mathcal Fu\Big).\label{fre}
\end{align}

where $\mathcal{F}$ is the Fourrier transform. Thanks to the Duhamel formula, solutions to \eqref{S} and \eqref{H} satisfy respectively the integral formula
\begin{align}
u(t)&:= e^{{\textnormal i}t\Delta}u_0- {\textnormal i}\int_0^t e^{{\textnormal i}(t-s)\Delta}\Big(|u|^{q-1}|x|^{-{b}}u\Big)\,ds;\label{duhamel}\\
w(t)&:= e^{{\textnormal i}t\Delta}u_0- {\textnormal i}\int_0^t e^{{\textnormal i}(t-s)\Delta}\Big((J_\alpha *|\cdot|^{-b}|w|^q)|w|^{q-2}|x|^{-b}w\Big)\,ds.\label{duhamel2}
\end{align}

Solutions of the problem \eqref{S} and \eqref{H}, formally satisfy the conservation laws

\begin{align}
M(u(t))&:= \|u(t)\|^2 = M(u_0),\quad M(w(t))= M(w_0)\tag{Mass};\\
E({u}(t))&:= \|\nabla u(t)\|^2+\frac2{1+q}\int_{\R^N} |u(t,x)|^{1+q}|x|^{-{b}}\,dx = E(u_0);\label{nrg}\tag{Energy}\\
E({w}(t))&:= \|\nabla w(t)\|^2+\frac1{q}\int_{\R^N}|w|^{q}(J_\alpha *|\cdot|^{-b}|w|^q)|x|^{-b}\,dx = E(w_0)\label{nrg2}\tag{Energy}.
\end{align}

If $u$ resolves the equation \eqref{S}, then also does the family ${u}_\kappa:=\kappa^\frac{2-{b}}{q-1}{u}(\kappa^{2}\cdot,\kappa \cdot)$, where $\kappa>0$. The critical exponent $s_c$ keeps invariant the following homogeneous Sobolev norm

\begin{align}
    \|{u}_\kappa(t)\|_{\dot H^\mu}=\kappa^{\mu-(\frac N2-\frac{2-b}{q-1})}\|{u}(\kappa^2 t)\|_{\dot H^\mu}:=\kappa^{\mu-s_c}\|{u}(\kappa^2 t)\|_{\dot H^\mu}.
\end{align}

The Shcr\"odinger problem \eqref{S} is said mass-critical when $s_c=0$ or $q=1+\frac{2(2-b)}{N}$. In the mass sub-critical regime $1<q< 1+\frac{2(2-b)}{N}$, it is known that for a datum with a finite mass, there is a global solution to \eqref{S} in the space $C([0,\infty),L^2)$, see \cite[Theorem 1.3]{cmg}. The equation \eqref{H} enjoys the scaling invariance

\begin{align}\label{scal}
{w}_\kappa:=\kappa^\frac{2-2b+\alpha}{2(q-1)}w(\kappa^{2}\cdot,\kappa\cdot),\quad\kappa>0.
\end{align}

The critical exponent $s_c'$ keeps invariant the following homogeneous Sobolev norm

\begin{align}
    \|w_\kappa(t)\|_{\dot H^\mu}=\kappa^{\mu-(\frac N2-\frac{2-2b+\alpha}{2(q-1)})}\|{w}(\kappa^2 t)\|_{\dot H^\mu}:=\kappa^{\mu-s_c'}\|{w}(\kappa^2 t)\|_{\dot H^\mu}.
\end{align}

The Hartree problem \eqref{H} is said mass-critical when $s_c'=0$ or $q=1+\frac{2-2b+\alpha}{N}$. In the mass sub-critical regime $2\leq q< 1+\frac{2-2b+\alpha}{N}$, it is known that for a datum with a finite mass, there is a global solution to \eqref{H} in the space $C([0,\infty),L^2)$, see \cite[Theorem 5.1]{tt}. 


For the reader convenience, we recall the so-called Hardy-Littlewood-Sobolev inequality \cite{el}.

\begin{lem}\label{hls}
Let $N\geq1$ and $0 <\alpha < N$.
\begin{enumerate}
\item[1.]
Let $s\geq1$ and $\frac1r=\frac1s+\frac\alpha N$. Then,
$$\|J_\alpha*g\|_s\leq C_{N,s,\alpha}\|g\|_{r},\quad\forall g\in L^r.$$
\item[2.]
Let $1<s,r,t<\infty$ be such that $\frac1r +\frac1s=\frac1t +\frac\alpha N$. Then,
$$\|(J_\alpha*g)f\|_t\leq C_{N,s,\alpha}\|f\|_{r}\|g\|_{s},\quad\forall (f,g)\in L^r\times L^s.$$
\end{enumerate}
\end{lem}

Eventually, let us give the dispersive estimate \cite{caz},

\begin{align}
    \|e^{it\Delta}\cdot\|_r\lesssim\frac1{t^{N(\frac12-\frac{1}r)}}\|\cdot\|_{r'},\quad\forall r\geq2.\label{dsp}
\end{align}

Now, we list the results proved in this work. 


\subsection{Main results}
We consider two cases.
~{\rm 
\begin{itemize}
    \item The inhomogeneous nonlinear Schr\"odinger problem \eqref{S}.
    
The first contribution of this note is the next non-scattering result.

\begin{thm}\label{t1}
Let $N\geq1$, $0<{b}<1$ and $1<q<1+\frac{2-2b}{N}$. Let a global solution of \eqref{S} denoted by ${u}\in C(\R_+, L^2)$. Assume that there is $u_+\in L^2$ such that

\begin{align}\label{t12}
\lim_{t\to\infty}\|u(t)-e^{it\Delta}u_+\|=0,
\end{align}

then, $u_+=0$.

\end{thm}

In view of the  results stated in the above theorem, some comments are in order.
~{\rm 
\begin{itemize}
\item 
The local existence of mass solutions to \eqref{S} was proved in \cite{gs,cmg,vdd}.
\item
By Theorem \ref{t2}, the assumption $q<1+\frac{2-2b}{N}$ is sharp.
\item
this work complements \cite{barab,sxa} to the inhomogeneous regime.
\end{itemize}}

The second contribution of this note is the next scattering result.

\begin{thm}\label{t2}
Let $N\geq1$, $0<{b}<1$ and $1+\frac{2-2b}{N}\leq q<1+\frac{2(2-b)}{N}$. Let $u_0\in \Sigma$ and the global solution ${u}\in C(\R_+, L^2)$ of \eqref{S}. Then, there is $u_+\in L^2$ such that

\begin{align}\label{t22}
\lim_{t\to\infty}\|u(t)-e^{it\Delta}u_+\|=0.
\end{align}

\end{thm}

In view of the  results stated in the above theorem, some comments are in order.
~{\rm 
\begin{itemize}
\item
We denote the weighted space $\Sigma:=\{f\in H^1, \,xf\in L^2\}$.
\item
This work complements \cite{ytky,kn} to the inhomogeneous case.
\item 
Using the standard virial type identity \cite[proposition 6.5.1]{caz},
\begin{align}
    \partial^2_{tt}\|xu(t)\|^2=8E(u_0)+4(N(q-1)+2b-4)\int_{\R^N}|u|^{1+q}|x|^{-b}\,dx\label{vrl},
\end{align}
the global energy solution with data in $\Sigma$ belongs to $C(\R,\Sigma)$.
\item 
The proof is given for $q<1+\frac{2(2-b)}{N}$. Indeed, in the mass super-critical regime, the scattering is known \cite{vdd}.
\end{itemize}}

\item The inhomogeneous nonlinear generalized Hartree problem \eqref{H}.

The third contribution of this note is the next non-scattering result.

\begin{thm}\label{t1'}
Let $N\geq2$, $b,\alpha$ satisfying \eqref{cnd} and $\max\{1,1+\frac{\alpha-2b}{N}\}<q<1+\frac{1+\alpha-2b}{N}$. Let a global solution of \eqref{H} denoted by ${w}\in C(\R_+, L^2)$ . Assume that there is $w_+\in L^2$ such that

\begin{align}\label{t12'}
\lim_{t\to\infty}\|w(t)-e^{it\Delta}w_+\|=0,
\end{align}

then, $w_+=0$.    
\end{thm}

In view of the  results stated in the above theorem, some comments are in order.
~{\rm 
\begin{itemize}
\item 
The local existence of solutions to \eqref{H} was proved in \cite{tt} under the weaker assumption $2-2b+\alpha>0$ rather than $1-2b+\alpha>0$.
\item
By Theorem \ref{t2'}, the assumption $q<1+\frac{1+\alpha-2b}{N}$ is sharp.
\item
The case $N=1$ follows under the extra assumption $q<{\frac32+\frac{\alpha-2b}{N}}$.
\item 
this work complements \cite{cgl} to the inhomogeneous case.
\end{itemize}}

The last contribution of this note is the next scattering result.

\begin{thm}\label{t2'}
Let $N\geq1$, $b,\alpha$ satisfying \eqref{cnd} and $1+\frac{1+\alpha-2b}{N}<q<1+\frac{2-2b+\alpha}{N}$. Let $w_0\in \Sigma$ and the global solution ${w}\in C(\R_+, L^2)$ of \eqref{H}. Then, there is $w_+\in L^2$ such that

\begin{align}\label{t22'}
\lim_{t\to\infty}\|w(t)-e^{it\Delta}w_+\|=0.
\end{align}

\end{thm}

In view of the  results stated in the above theorem, some comments are in order.
~{\rm 
\begin{itemize}
\item
this work complements \cite{cgl} to the inhomogeneous case and $q\neq2$.
\item 
Using the standard virial type identity \cite[proposition 2.7]{tt},
\begin{align}
    \partial^2_{tt}\|xw(t)\|^2=8\|\nabla w\|^2+\frac{4(N(q-1)-\alpha-2b)}{q}\int_{\R^N}(J_\alpha*|\cdot|^{-b}|w|^q)|x|^{-b}|w|^q\,dx \label{vrl'},
\end{align}
the global energy solution with data in $\Sigma$ belongs to $C(\R,\Sigma)$.
\item 
The proof is given for $q<1+\frac{2-2b+\alpha}{N}$. Indeed, in the mass super-critical regime, the scattering is known \cite{subm}.
\end{itemize}}
\end{itemize}}

In the rest of this note, we prove the main results.

\section{Non scattering}\label{sec3}

In this section, we prove Theorem \ref{t1} and Theorem \ref{t1'}.

\subsection{Non scattering of the inhomogeneous nonlinear Schr\"odinger equation}\label{sec41}

In this sub-section, we prove Theorem \ref{t1}. Assume that \eqref{t12} holds and denote by $v:=e^{i\cdot\Delta}u_+$. Then, for any $t\geq0$, yields

\begin{align}
    \|u_0\|=\|v(t)\|=\|u_+\|.\label{01}
\end{align}

We define the real function

\begin{align}
    \Upsilon:=\Upsilon_{u,v}: t\mapsto\Im\int_{\R^N}\bar uv\,dx.\label{02}
\end{align}

Let us denote also 

\begin{align}
\Theta(f,g):=-|x|^{-b}|f|^{q-1}g,\quad \Theta(f):=\Theta(f,f).\label{04}    
\end{align}

A direct computation via \eqref{S} gives

\begin{align}
    \Upsilon'(t)
    &=\Im\int_{\R^N}\big(\bar u\partial_tv-\partial_tu\bar v\big)\,dx\nonumber\\
        &=\Im\int_{\R^N}\big(i\Delta v\bar u-i(\Delta u+\Theta(u))\bar v\big)\,dx\nonumber\\
        &=\Re\int_{\R^N}\big(\Delta v\bar u-(\Delta u+\Theta(u))\bar v\big)\,dx\label{03}
        \end{align}

So, \eqref{03} implies that 

   \begin{align}
    \Upsilon'(t)     
        &=\Re\int_{\R^N}|x|^{-b}|u|^{q-1}u\bar v\,dx\nonumber\\
                &=\Re\int_{\R^N}|x|^{-b}\big(|u|^{q-1}u+|v|^{q-1}v-|v|^{q-1}v\big)\bar v\,dx\nonumber\\
                &=\int_{\R^N}|v|^{1+q}|x|^{-b}\,dx+\Re\int_{\R^N}|x|^{-b}\big(|u|^{q-1}u-|v|^{q-1}v\big)\bar v\,dx\label{05}.
\end{align}

Letting $\alpha>0$, we get by H\"older estimate

\begin{align}
   \|v(t)\|^2_{L^2(\Omega_{\alpha t})}
    &=\||x|^\frac{b}{1+q}|x|^\frac{-b}{1+q}v(t)\|_{L^2(\Omega_{\alpha t})}^2\nonumber\\
      &\leq\||x|^\frac{b}{1+q}\|_{L^\frac{2(1+q)}{q-1}(\Omega_{\alpha t})}^2\||x|^\frac{-b}{1+q}v(t)\|_{L^{1+q}(\Omega_{\alpha t})}^2\nonumber\\
      &\lesssim (\alpha t)^\frac{N(q-1)+2b}{1+q}\||x|^\frac{-b}{1+q}v(t)\|_{L^{1+q}(\Omega_{\alpha t})}^2\nonumber\\
      &\lesssim (\alpha t)^\frac{N(q-1)+2b}{1+q}\big(\int_{\R^N}|x|^{-b}|v(t)|^{1+q}\,dx\big)^\frac{2}{1+q}\label{06}.
\end{align}

Now, we write

\begin{align}
    v(t)
    &=\mathcal{F}^{-1}(e^{it|\cdot|^2})*u_+\nonumber\\
    &=ct^{-\frac{N}{2}}e^{\frac{i|\cdot|^2}{4t}}*u_+.\label{08}
\end{align}

So, \eqref{08} implies that

\begin{align}
    \|v(t)\|^2_{L^2(\Omega_{\alpha t})}
    &=t^{-N}\int_{\Omega_{\alpha t}}|e^{\frac{i|\cdot|^2}{4t}}*u_+|^2\,dx\nonumber\\
    &=t^{-N}\int_{\Omega_{\alpha t}}\big|\int_{\R^N}e^{\frac{i|x-y|^2}{4t}}u_+(y)\,dy\big|^2\,dx\nonumber\\
    &=t^{-N}\int_{\Omega_{\alpha t}}\big|\int_{\R^N}e^{-i\frac{x\cdot y}{2t}}e^{\frac{i|y|^2}{4t}}u_+(y)\,dy\big|^2\,dx\label{09}.
\end{align}

Thanks to the variable change $z=-\frac1{2t}x$, the equality \eqref{09} gives for some $c>0$,

\begin{align}
    \|v(t)\|^2_{L^2(\Omega_{\alpha t})}
    &=c\int_{\Omega_{\alpha/2 }}\big|\int_{\R^N}e^{i{z\cdot y}}e^{\frac{i|y|^2}{t}}u_+(y)\,dy\big|^2\,dz\nonumber\\
        &=c\int_{\Omega_{\alpha/2}}\big|\mathcal{F}^{-1}\big(e^{\frac{i|\cdot|^2}{t}}u_+\big)\,\big|^2\,dz\label{010}.
\end{align}

Moreover, by the fact that $\mathcal F$ is an isomorphism of $L^2$, it follows that by \eqref{010},

\begin{align}
    \lim_{t\to\infty}\|v(t)\|^2_{L^2(\Omega_{\alpha t})}
        &=c\int_{\Omega_{\alpha/2}}\big|\big(\mathcal{F}^{-1}(u_+)\big)(z)\,\big|^2\,dz\label{011}.
\end{align}

So, \eqref{011} implies that for $\alpha\gg1$, holds 

\begin{align}
    \lim_{t\to\infty}\|v(t)\|^2_{L^2(\Omega_{\alpha t})}
        &\gtrsim\int_{\R^N}\big|\big(\mathcal{F}^{-1}(u_+)\big)(z)\,\big|^2\,dz=\|u_+\|^2\label{012}.
\end{align}

Then, \eqref{012} implies that for $t>T\gg1$, we have

\begin{align}
    \|v(t)\|^2_{L^2(\Omega_{\alpha t})}
        &\gtrsim\|u_+\|^2\gtrsim 1\label{013}.
\end{align}

Now, \eqref{06} and \eqref{013} give

\begin{align}
   1&\lesssim \|v(t)\|^2_{L^2(\Omega_{\alpha t})}\nonumber\\
      &\lesssim (\alpha t)^\frac{N(q-1)+2b}{1+q}\big(\int_{\R^N}|v|^{1+q}|x|^{-b}\,dx\big)^\frac{2}{1+q}\label{014}.
\end{align}

Finally, we rewrite \eqref{014} as

\begin{align}
   \int_{\R^N}|v|^{1+q}|x|^{-b}\,dx&\gtrsim  (\alpha t)^{-\frac{N(q-1)+2b}{2}},\quad\forall t>T\label{016}.
\end{align}

Let us write

   \begin{align}
      \big|\int_{\R^N}|x|^{-b}\big(|u|^{q-1}u-|v|^{q-1}v\big)\bar v\,dx\big|
      &\lesssim\int_{\Omega}|x|^{-b}\big(|u|^{q-1}+|v|^{q-1}\big)|u-v||v|\,dx\nonumber\\
      &+\int_{\Omega^c}|x|^{-b}\big(|u|^{q-1}+|v|^{q-1}\big)|u-v||v|\,dx\nonumber\\
      &:=(I)+(II).\label{015}
\end{align}

Using H\"older estimate, we write

   \begin{align}
         (I)
      &=\int_{\Omega}|x|^{-b}\big(|u|^{q-1}+|v|^{q-1}\big)|u-v||v|\,dx\nonumber\\
&\leq \||x|^{-b}\|_{L^a(\Omega)}\big(\|u\|^{q-1}+\|v\|^{q-1}\big)\|u-v\|\|v\|_\gamma\label{017}.
\end{align}

Here,

\begin{align}
    \frac bN<\frac1a=1-\frac q2-\frac1\gamma.\label{018}
\end{align}

Moreover, \eqref{018} reads

\begin{align}
     0<\frac1\gamma<1-\frac q2-\frac bN.\label{019}
\end{align}

The inequality \eqref{019} is satisfied if

\begin{align}
    q<1+\frac2N-\frac{2b}N.\label{020}
\end{align}

Moreover, by H\"older estimate

\begin{align}
         (II)
      &=\int_{\Omega^c}|x|^{-b}\big(|u|^{q-1}+|v|^{q-1}\big)|u-v||v|\,dx\nonumber\\
&\leq \||x|^{-b}\|_{L^d(\Omega^c)}\big(\|u\|^{q-1}+\|v\|^{q-1}\big)\|u-v\|\|v\|_\beta\label{21}.
\end{align}

Here,

\begin{align}
    \frac bN>\frac1d=1-\frac q2-\frac1\beta.\label{22}
\end{align}

Moreover, \eqref{22} reads

\begin{align}
    \frac{1}{2}\geq\frac1\beta>-\frac q2-\frac bN.\label{23}
\end{align}

The inequality \eqref{23} is rewritten as

\begin{align}
    q>1-\frac{2b}N.\label{24}
\end{align}

Then, \eqref{020} and \eqref{24} read

\begin{align}
   1-\frac{2b}N<  q<1+\frac2N-\frac{2b}N.\label{25}
\end{align}

Now, under the assumption \eqref{25}, we get by \eqref{015}, \eqref{017} and \eqref{21}, via \eqref{dsp},

\begin{align}
      \int_{\R^N}|x|^{-b}\big(|u|^{q-1}u-|v|^{q-1}v\big)\bar v\,dx
      &\lesssim \big(\|v\|_\gamma+\|v\|_\beta\big)\|u-v\|\nonumber\\
      &\lesssim \big(t^{-N(\frac12-\frac{1}{\gamma})}+t^{-N(\frac12-\frac{1}{\beta})}\big)\|u-v\|.\label{26}
\end{align}

Then, by \eqref{26} via \eqref{05} and \eqref{016}, we get for all $t>T$,

\begin{align}
    \Upsilon'(t)
    &=\int_{\R^N}|v|^{1+q}|x|^{-b}\,dx+\Re\int_{\R^N}|x|^{-b}\big(|u|^{q-1}u-|v|^{q-1}v\big)\bar v\,dx\nonumber\\
    &\gtrsim (\alpha t)^{-\frac{N(q-1)+2b}{2}}-\big(t^{-N(\frac12-\frac{1}{\gamma})}+t^{-N(\frac12-\frac{1}{\beta})}\big)\|u-v\|\nonumber\\
    &\gtrsim (\alpha t)^{-\frac{N(q-1)+2b}{2}}\label{27}.
\end{align}

Indeed, we have by \eqref{019},

\begin{align}
    N(\frac12-\frac1\gamma)
    &>N(\frac12-1+\frac q2+\frac bN)\nonumber\\
    &\geq\frac{N(q-1)+2b}{2},\label{127}
\end{align}

and in \eqref{23}, taking $\beta\gg1$, we get for ${N\geq2}$,

\begin{align}
    N(\frac12-\frac1\beta)
    &\simeq \frac N2\nonumber\\
    &\geq\frac{N(q-1)+2b}{2}.\label{127'}
\end{align}

Integrating \eqref{27}, we get 

\begin{align}
\sup_{t>T}\Upsilon(t)
        &\gtrsim (\alpha T)^{1-\frac{N(q-1)+2b}{2}}\nonumber\\
        &=\infty,\label{28}
\end{align}

where, we used \eqref{25}, namely

\begin{align}
    \frac{N(q-1)+2b}{2}<1.\label{29}
\end{align}

The contradiction \eqref{28} ends the proof.

\subsection{Non scattering of the inhomogeneous nonlinear generalized Hartree equation}\label{sec42}

In this sub-section, we prove Theorem \ref{t1'}.  Assume that \eqref{t12'} holds and denote by $v:=e^{i\cdot\Delta}w_+$. Then, for any $t\geq0$, yields

\begin{align}
    \|w_0\|=\|v(t)\|=\|w_+\|.\label{01'}
\end{align}

We define the real function

\begin{align}
    \Upsilon:=\Upsilon_{w,v}: t\mapsto\Im\int_{\R^N}\bar wv\,dx.\label{02'}
\end{align}

Let us denote also 

\begin{align}
\Theta:=\Theta(f,g):=-(J_\alpha*|\cdot|^{-b}|f|^q)|x|^{-b}|f|^{q-2}g,\quad \Theta(f):=\Theta(f,f).\label{04'}    
\end{align}

A direct computation via \eqref{H} gives

\begin{align}
    \Upsilon'(t)
    &=\Im\int_{\R^N}\big(\bar w\partial_tv-\partial_tw\bar v\big)\,dx\nonumber\\
        &=\Im\int_{\R^N}\big(i\Delta v\bar w-i(\Delta w+\Theta(w))\bar v\big)\,dx\nonumber\\
        &=\Re\int_{\R^N}\big(\Delta v\bar w-(\Delta w+\Theta(w))\bar v\big)\,dx\label{03'}
        \end{align}

So, \eqref{03'} implies that 

   \begin{align}
    \Upsilon'(t)     
        &=\Re\int_{\R^N}(J_\alpha*|\cdot|^{-b}|w|^q)|w|^{q-2}|x|^{-b}w\bar v\,dx\nonumber\\
                &=\Re\int_{\R^N}(J_\alpha*|\cdot|^{-b}|w|^q)|x|^{-b}\big(|w|^{q-2}w+|v|^{q-2}v-|v|^{q-2}v\big)\bar v\,dx\nonumber.
                \end{align}

                Thus,

                \begin{align}
    \Upsilon'(t) 
                &=\int_{\R^N}(J_\alpha*|\cdot|^{-b}|w|^q)|x|^{-b}|v|^q\,dx+\Re\int_{\R^N}(J_\alpha*|\cdot|^{-b}|w|^q)|x|^{-b}\big(|w|^{q-2}w-|v|^{q-2}v\big)\bar v\,dx\nonumber\\
                &=\int_{\R^N}(J_\alpha*|\cdot|^{-b}|v|^q)|x|^{-b}|v|^q\,dx+\int_{\R^N}(J_\alpha*|\cdot|^{-b}[|w|^q-|v|^q])|x|^{-b}|v|^q\,dx\nonumber\\
                &+\Re\int_{\R^N}(J_\alpha*|\cdot|^{-b}|w|^q)|x|^{-b}\big(|w|^{q-2}w-|v|^{q-2}v\big)\bar v\,dx\nonumber\\
                &:=(A)+(B)+(C).\label{05'}
\end{align}

Taking $\gamma>0$, we write for $|x|<\gamma t$,

\begin{align}
(J_\alpha*|\cdot|^{-b}|v|^q)(x)
&=C_{N,\alpha}\int_{\R^N}\frac{|y|^{-b}|v(y)|^q}{|x-y|^{N-\alpha}}\,dy\nonumber\\
&\gtrsim\int_{\Omega_{\gamma t}}\frac{|y|^{-b}|v(y)|^q}{|x-y|^{N-\alpha}}\,dy\nonumber\\
&\gtrsim t^{-(N-\alpha)-b}\int_{\Omega_{\gamma t}}{|v(y)|^q}\,dy\label{051}.
\end{align}

So, \eqref{05'} gives via \eqref{051},

\begin{align}
    (A)
    &\geq \int_{\Omega_{\gamma t}}(J_\alpha*|\cdot|^{-b}|v|^q)|x|^{-b}|v|^q\,dx\nonumber\\
        &\gtrsim t^{-(N-\alpha)-b}\int_{\Omega_{\gamma t}}\big(\int_{\Omega_{\gamma t}}{|v(y)|^q}\,dy\big)|x|^{-b}|v(x)|^q\,dx\nonumber\\
        &\gtrsim t^{-(N-\alpha)-2b}\big(\int_{\Omega_{\gamma t}}|v(x)|^q\,dx\big)^2\label{052}.
\end{align}

Moreover, by H\"older estimate

\begin{align}
    \|v\|_{L^2(\Omega_{\gamma t})}
    &\leq \|v\|_{L^q(\Omega_{\gamma t})}|\Omega_{\gamma t}|^\frac{q-2}{2q}\nonumber\\
    &\lesssim \|v\|_{L^q(\Omega_{\gamma t})}(\gamma t)^\frac{N(q-2)}{2q}\label{053}.
\end{align}

So, by \eqref{052} and \eqref{053} via \eqref{013}, it follows that

\begin{align}
    (A)
        &\gtrsim t^{-(N-\alpha)-2b}\big((\gamma t)^{-\frac{N(q-2)}{2q}}\|v\|_{L^2(\Omega_{\gamma t})}\big)^{2q}\nonumber\\
        &\gtrsim \gamma ^{-N(q-2)}t^{-(N(q-1)-\alpha)-2b}\|v\|_{L^2(\Omega_{\gamma t})}^{2q}\nonumber\\
        &\gtrsim \gamma ^{-N(q-2)}t^{-(N(q-1)-\alpha)-2b}\label{054}.
\end{align}

Now, we estimate $(C)$. By Lemma \ref{hls}, we write

\begin{align}
    (C_1)
                &:=\big|\int_{\Omega}(J_\alpha*|\cdot|^{-b}|w|^q)|x|^{-b}\big(|w|^{q-2}w-|v|^{q-2}v\big)\bar v\,dx\big|\nonumber\\
                &\lesssim \int_{\Omega}(J_\alpha*|\cdot|^{-b}|w|^q)|x|^{-b}\big(|w|^{q-2}+|v|^{q-2}\big)|v| |v-w|\,dx\nonumber\\
                &\lesssim \big(\||x|^{-b}\|_{L^{e_1}(\Omega)}\||x|^{-b}\|_{L^{f_1}(\Omega)}+\||x|^{-b}\|_{L^{f_2}(\Omega^c)}\||x|^{-b}\|_{L^{e_2}(\Omega)}\big)\|w\|^q\big(\|w\|^{q-2}+\|v\|^{q-2}\big)\|v\|_r \|v-w\|\nonumber\\
                &\lesssim \big(\|w\|^{q-2}+\|v\|^{q-2}\big)\|w\|^q\|v\|_r \|v-w\|\label{055}.
\end{align}

Here

\begin{equation}
\left\{
\begin{array}{ll}
1+\frac{\alpha}N=\frac1{e_1}+\frac1{f_1}+\frac{2q-1}2+\frac1{r}=\frac1{e_2}+\frac1{f_2}+\frac{2q-1}2+\frac1{r} ;\\
\min\{\frac N{e_1},\frac N{e_2},\frac N{f_1}\}>{b}> \frac N{f_2};\\
r>2.
\label{m-1'}
\end{array}
\right.
\end{equation}

So, we need

\begin{align}
    \frac32+\frac{\alpha}N-q-\frac{1}{r}=\frac1{e_1}+\frac1{f_1}>\frac{2b}{N}.\label{056}
\end{align}

Then, \eqref{056} reads

\begin{align}
    \frac{1}{r}<\frac32+\frac{\alpha}N-q-\frac{2b}{N}.\label{057}
\end{align}

Thus, \eqref{055} holds, since \eqref{057} is possible because for $N\geq2$,

\begin{align}
{q<1+\frac{1+\alpha-2b}{N}\leq\frac32+\frac{\alpha-2b}{N}}.    
\end{align}

Also by Lemma \ref{hls}, we write

\begin{align}
    (C_2)
                &:=\big|\int_{\Omega^c}(J_\alpha*|\cdot|^{-b}|w|^q)|x|^{-b}\big(|w|^{q-2}w-|v|^{q-2}v\big)\bar v\,dx\big|\nonumber\\
                &\lesssim \int_{\Omega^c}(J_\alpha*|\cdot|^{-b}|w|^q)|x|^{-b}\big(|w|^{q-2}+|v|^{q-2}\big)|v| |v-w|\,dx\nonumber\\
                &\lesssim \big(\||x|^{-b}\|_{L^{g_1}(\Omega)}\||x|^{-b}\|_{L^{h_1}(\Omega^c)}+\||x|^{-b}\|_{L^{g_2}(\Omega^c)}\||x|^{-b}\|_{L^{h_2}(\Omega^c)}\big)\|w\|^q\big(\|w\|^{q-2}+\|v\|^{q-2}\big)\|v\|_{r_2}\|v-w\|\nonumber\\
                &\lesssim \big(\|w\|^{q-2}+\|v\|^{q-2}\big)\|w\|^q\|v\|_{r_2} \|v-w\|\label{0552}.
\end{align}

Here

\begin{equation}
\left\{
\begin{array}{ll}
1+\frac{\alpha}N=\frac1{g_1}+\frac1{h_1}+\frac{2q-1}2+\frac1{r_2}=\frac1{g_2}+\frac1{h_2}+\frac{2q-1}2+\frac1{r_2} ;\\
\min\{\frac N{h_1},\frac N{h_2},\frac N{g_2}\}<{b}< \frac N{g_1};\\
r_2>2.
\label{m-1''}
\end{array}
\right.
\end{equation}

So, we need

\begin{align}
    \frac32+\frac{\alpha}N-q-\frac1{r_2}=\frac1{h_1}+\frac1{g_1}<\frac{2b}{N}.\label{0562}
\end{align}

Then, \eqref{0562} reads

\begin{align}
\frac32+\frac{\alpha}N-q-\frac{2b}{N}<\frac{1}{r_2}<\frac12.\label{0572}
\end{align}

Thus, \eqref{0552} holds, since \eqref{0572} is possible for 

\begin{align}
{q>1+\frac{\alpha-2b}{N}}.    
\end{align}

Since the term $(B)$ can be estimated like $(C)$, we write by \eqref{055} and \eqref{0552} via \eqref{dsp}, 

\begin{align}
    (B)+(C)\lesssim t^{-N(\frac12-\frac1r)}+t^{-N(\frac12-\frac1{r_2})}.\label{058}
\end{align}

Thus, \eqref{05'} via \eqref{054} and \eqref{058}, we write for $t\gg1\gg\varepsilon>0$,

\begin{align}
    \Upsilon'(t) 
                &\gtrsim t^{-(N(q-1)-\alpha)-2b}-\big(t^{-N(\frac12-\frac1r)}+t^{-N(\frac12-\frac1{r_2})}\big)\nonumber\\
                 &\gtrsim t^{-(N(q-1)-\alpha)-2b+\varepsilon}.\label{059}
\end{align}

Indeed, by \eqref{057}, we have

\begin{align}
    N(\frac12-\frac1r)
    &>N(\frac12-\frac32-\frac\alpha N+q+\frac{2b}{N})\nonumber\\
    &>N(q-1)-\alpha+2b.\label{1059}
\end{align}

Also, by \eqref{0572}, taking $\frac1{r_2}$ close to $\frac32+\frac{\alpha}N-q-\frac{2b}{N}$, we have 

\begin{align}
N(\frac12-\frac1{r_2})
   &\simeq N\big( \frac12-(\frac32+\frac{\alpha}N-q-\frac{2b}{N})\big)\nonumber\\
    &=N(q-1)-{\alpha}+2b.
\end{align}

Now, since 

\begin{align}
    {q<1+\frac{1+\alpha-2b}{N}},
\end{align}

an integration of \eqref{059} gives

\begin{align}
    \lim_{t\to\infty}\Upsilon(t)=\infty.\label{060}
\end{align}

Since $\Upsilon(t)\leq \|w_0\|^2$, the limit \eqref{060} gives a contradiction, which ends the proof. 

\section{Scattering}\label{sec4}

In this section, we prove Theorem \ref{t2} and Theorem \ref{t2'}.

\subsection{Scattering of the inhomogeneous nonlinear Schr\"odinger equation}\label{sec31}

In this sub-section, we prove Theorem \ref{t2}. Take the transformation

\begin{align}
    v(t,x)&:=(it)^{-\frac N2}e^{i\frac{|x|^2}{4t}}\bar u(\frac{1}{t},\frac{x}{t}).\label{30}
\end{align}

We compute the time derivative

\begin{align}
    \partial_tv(t,x)
    &=\partial_t\big((it)^{-\frac N2}e^{i\frac{|x|^2}{4t}}\bar u(\frac{1}{t},\frac{x}{t})\big)\nonumber\\
    &=(it)^{-\frac N2}e^{i\frac{|x|^2}{4t}}\big([-\frac N{2t}-i\frac{|x|^2}{4t^2}]\bar u(\frac{1}{t},\frac{x}{t})-\frac{1}{t^2}{\partial_t}\bar u(\frac{1}{t},\frac xt)-\frac{1}{t^2}x\cdot\nabla \bar u(\frac{1}{t},\frac{x}{t})\big)\label{31}.
\end{align}

Moreover, the Laplacian reads

\begin{align}
    \Delta v(t,x)
    &=(it)^{-\frac N2}\Delta\big(e^{i\frac{|x|^2}{4t}}\bar u(\frac{1}{t},\frac{x}{t})\big)\nonumber\\
    &=(it)^{-\frac N2}e^{i\frac{|x|^2}{4t}}\big(\frac{1}{t^2}\Delta\bar u(\frac{1}{t},\frac{x}{t})+\frac{i}{t^2}x\cdot\nabla\bar u(\frac{1}{t},\frac{x}{t})+[\frac{Ni}{2t}-\frac{|x|^2}{4t^2}]\bar u(\frac{1}{t},\frac{x}{t})  \big)\label{32}.
\end{align}

Taking account of \eqref{31} and \eqref{32}, we get

\begin{align}
    i\partial_tv(t,x)+\Delta v
        &=(it)^{-\frac N2}e^{i\frac{|x|^2}{4t}}\Big(i\big([-\frac N{2t}-i\frac{|x|^2}{4t^2}]\bar u(\frac{1}{t},\frac{x}{t})-\frac{1}{t^2}{\partial_t}\bar u(\frac{1}{t},\frac xt)-\frac{1}{t^2}x\cdot\nabla \bar u(\frac{1}{t},\frac{x}{t})\big)\nonumber\\
    &+\big(\frac{1}{t^2}\Delta\bar u(\frac{1}{t},\frac{x}{t})+\frac{i}{t^2}x\cdot\nabla\bar u(\frac{1}{t},\frac{x}{t})+[\frac{Ni}{2t}-\frac{|x|^2}{4t^2}]\bar u(\frac{1}{t},\frac{x}{t})  \big)\Big)\nonumber\\
    &=-(it)^{-\frac N2-2}e^{i\frac{|x|^2}{4t}}\Big(-i{\partial_t}\bar u(\frac{1}{t},\frac xt)+\Delta\bar u(\frac{1}{t},\frac{x}{t})  \Big)\label{33}.
\end{align}

So, \eqref{33} via \eqref{S} and \eqref{30} gives

\begin{align}
    i\partial_tv(t,x)+\Delta v
    &=-(it)^{-\frac N2-2}e^{i\frac{|x|^2}{4t}}\big(\frac{|x|}t\big)^{-b}|u(\frac{1}{t},\frac{x}{t})|^{q-1}\bar u(\frac{1}{t},\frac{x}{t})  \nonumber\\
    &=-(it)^{-2}\big(\frac{|x|}t\big)^{-b}|t^\frac{N}{2}v|^{q-1}v  \nonumber\\
    &=t^{\frac{N(q-1)}{2}+b-2}|x|^{-b}|v|^{q-1}v  \label{34}.
\end{align}

We test \eqref{34} with $\frac1{t^{\frac{N(q-1)}{2}+b-2}}\partial_t\bar v$, as follows

\begin{align}
    \Re\Big(\frac1{t^{\frac{N(q-1)}{2}+b-2}}\partial_t\bar v\big(i\partial_tv(t,x)+\Delta v\big)\Big)
    &=\Re\big(\partial_t\bar v|x|^{-b}|v|^{q-1}v\big) ,
\end{align}

which reads via the fact that $\frac{N(q-1)}{2}+b-2<0$,

\begin{align}
    \partial_t\big(\frac1{t^{\frac{N(q-1)}{2}+b-2}}\|\nabla v\|^2\big)
    &\geq\frac1{t^{\frac{N(q-1)}{2}+b-2}}\partial_t\|\nabla v\|^2\nonumber\\
    &=-\frac2{1+q}\partial_t\||x|^{-\frac b{1+q}}v\|_{1+q}^{1+q}.\label{35}
\end{align}

Integrating in time \eqref{35}, it follows that for $0<\tau\leq t<\infty$, 

\begin{align}
   \frac1{t^{\frac{N(q-1)}{2}+b-2}}\|\nabla v(t)\|^2-\frac1{\tau^{\frac{N(q-1)}{2}+b-2}}\|\nabla v(\tau)\|^2
    &\geq-\frac2{1+q}\big(\||x|^{-\frac b{1+q}}v(t)\|_{1+q}^{1+q}-\||x|^{-\frac b{1+q}}v(\tau)\|_{1+q}^{1+q}\big)\,\nonumber
\end{align}

which reads for any $0<\tau\leq t<\infty$,  

\begin{align}
   \frac1{t^{\frac{N(q-1)}{2}+b-2}}\|\nabla v(t)\|^2+\frac2{1+q}\||x|^{-\frac b{1+q}}v(t)\|_{1+q}^{1+q}
    &\geq\frac1{\tau^{\frac{N(q-1)}{2}+b-2}}\|\nabla v(\tau)\|^2+\||x|^{-\frac b{1+q}}v(\tau)\|_{1+q}^{1+q}\label{36}.
\end{align}

So, by \eqref{36} via the mass conservation law, we get

\begin{align}
  \max\big\{\sup_{0<t\leq1}\frac1{t^{\frac{N(q-1)}{2}+b-2}}\|\nabla v(t)\|^2,\sup_{0<t\leq1}\|v(t)\|,\sup_{0<t\leq1}\||x|^{-\frac b{1+q}}v(t)\|_{1+q}\big\}<\infty.\label{37}
\end{align}

Take a test function $\zeta\in H^1$ and compute

\begin{align}
   \int_{\R^N} (v(t)-v(\tau))\zeta\,dx
   &=\int_\tau^t\int_{\R^N}\partial_sv(s) \zeta\,dx\,ds\nonumber\\
      &=i\int_\tau^t\int_{\R^N}\big(\Delta v-s^{\frac{N(q-1)}{2}+b-2}|x|^{-b}|v|^{q-1}v \big)\zeta\,dx\,ds\nonumber\\
            &=-i\int_\tau^t\int_{\R^N}\nabla v\cdot\nabla\zeta\,dx\,ds-i\int_\tau^t\int_{\R^N}s^{\frac{N(q-1)}{2}+b-2}|x|^{-b}|v|^{q-1}v\zeta\,dx\,ds\label{38}.
\end{align}

Moreover, \eqref{38} via the fact that $\frac{N(q-1)}{2}+b-2>-1$, we get, when $t\to0$,

\begin{align}
    v(t)\rightharpoonup v_+,\quad\mbox{in}\quad L^2.\label{39}
\end{align}

Now, we pick $\zeta:=v(t)$ and rewrite \eqref{38} as follows

\begin{align}
   \big|\int_{\R^N} (v(t)-v(\tau))v(t)\,dx\big|
            &= \big|-i\int_\tau^t\int_{\R^N}\nabla v(s)\cdot\nabla v(t)\,dx\,ds\nonumber\\
            &-i\int_\tau^t\int_{\R^N}s^{\frac{N(q-1)}{2}+b-2}|x|^{-b}|v(s)|^{q-1}v(s)v(t)\,dx\,ds \big|\nonumber\\
            &\leq\int_\tau^t\|\nabla v(s)\|\|\nabla v(t)\|\,ds+\int_\tau^ts^{\frac{N(q-1)}{2}+b-2}\||x|^{-\frac{b}{1+q}}v(s)\|_{1+q}^{q}\||x|^{-\frac{b}{1+q}}v(t)\|_{1+q}\,ds \label{40}
\end{align}

Using \eqref{37} in \eqref{40}, we get

\begin{align}
   \big|\int_{\R^N} (v(t)-v(\tau))v(t)\,dx\big|
            &\lesssim t^{\frac{N(q-1)+2b}{4}-1}\int_\tau^ts^{\frac{N(q-1)+2b}{4}-1}\,ds+\int_\tau^ts^{\frac{N(q-1)}{2}+b-2}\,ds \nonumber\\
             &\lesssim t^{\frac{N(q-1)+2b}{4}-1}\big(t^{\frac{N(q-1)+2b}{4}}-\tau^{\frac{N(q-1)+2b}{4}}\big)+t^{\frac{N(q-1)+2b}{2}-1}-\tau^{\frac{N(q-1)+2b}{2}-1} .\label{41}
\end{align}

Hence, letting $\tau\to0$ in \eqref{41} and using \eqref{39}, we write

\begin{align}
   \big|\int_{\R^N} (v(t)-v(0))v(t)\,dx\big|
             &\lesssim t^{\frac{N(q-1)+2b}{2}-1}.\label{42}
\end{align}

Now, the triangular estimate via \eqref{42} gives

\begin{align}
\|v(t)-v(0)\|^2
   &\leq\big|\int_{\R^N} (v(t)-v(0))v(t)\,dx\big|+\big|\int_{\R^N} (v(t)-v(0))v(0)\,dx\big|\nonumber\\
             &\lesssim t^{\frac{N(q-1)+2b}{2}-1}+\big|\int_{\R^N} (v(t)-v(0))v(0)\,dx\big|\nonumber\\
             &\to0,\quad\mbox{as}\quad t\to0.\label{43}
\end{align}

So, \eqref{43} via \eqref{30}, implies that for $u_+(x):=e^{it\Delta}\big((it)^{-\frac N2}e^{i\frac{|x|^2}{4t}}\bar v_0(\frac xt)\big)$,

\begin{align}
    \|u(t)-e^{-it\Delta}u_+\|
    &=\|(it)^{-\frac N2}e^{i\frac{|x|^2}{4t}}\bar v(\frac{1}{t},\frac{x}{t})-e^{-it\Delta}u_+\|\nonumber\\
        &=\|(it)^{-\frac N2}e^{i\frac{|x|^2}{4t}}\bar v(\frac{1}{t},\frac{x}{t})-(it)^{-\frac N2}e^{i\frac{|x|^2}{4t}}\bar v_0(\frac xt)\|\nonumber\\
    &=\| v(\frac{1}{t},\cdot)-v(0)\|\nonumber\\
    &\to0,\quad\mbox{as}\quad t\to\infty.\label{44}
\end{align}

The scattering is proved thanks to \eqref{44}. Theorem \ref{t2} is established.

\subsection{Scattering of the inhomogeneous nonlinear generalized Hartree equation}\label{sec32}

In this sub-section, we prove Theorem \ref{t2'}. Take the transformation

\begin{align}
    v(t,x)&:=(it)^{-\frac N2}e^{i\frac{|x|^2}{4t}}\bar w(\frac{1}{t},\frac{x}{t}).\label{302}
\end{align}

By \eqref{33} via \eqref{H} and \eqref{30}, yields

\begin{align}
    i\partial_tv(t,x)+\Delta v
    &=-(it)^{-\frac N2-2}e^{i\frac{|x|^2}{4t}}\Big(-i{\partial_t}\bar w(\frac{1}{t},\frac xt)+\Delta\bar w(\frac{1}{t},\frac{x}{t})  \Big)\nonumber\\
    &=-(it)^{-\frac N2-2}e^{i\frac{|x|^2}{4t}}\big(\frac{|x|}t\big)^{-b}|w(\frac{1}{t},\frac{x}{t})|^{q-2}(J_\alpha*|\cdot|^{-b}|w|^q)(\frac1t,\frac xt)\bar w(\frac{1}{t},\frac{x}{t})  \nonumber\\
    &=-(it)^{-\frac N2-2}e^{i\frac{|x|^2}{4t}}\big(\frac{|x|}t\big)^{-b}|w(\frac{1}{t},\frac{x}{t})|^{q-2}t^{b-\alpha}(J_\alpha*|\cdot|^{-b}|w(\frac1t,\frac xt)|^q)\bar w(\frac{1}{t},\frac{x}{t}) \label{342'}.
    \end{align}

Thus, by \eqref{342'}, we write

    \begin{align}
    i\partial_tv(t,x)+\Delta v
    &=-(it)^{-2}\big(\frac{|x|}t\big)^{-b}|t^\frac{N}{2}v|^{q-2}t^{b-\alpha}(J_\alpha*|\cdot|^{-b}|t^\frac{N}{2}v|^q)v  \nonumber\\
    &=t^{N(q-1)+2b-\alpha-2}|x|^{-b}|v|^{q-2}(J_\alpha*|\cdot|^{-b}|v|^q)v  \label{342}.
\end{align}

We test \eqref{342} with $\frac1{t^{N(q-1)+2b-\alpha-2}}\partial_t\bar v$, as follows

\begin{align}
    \Re\Big(\frac1{t^{N(q-1)+2b-\alpha-2}}\big(i\partial_tv(t,x)+\Delta v\big)\partial_t\bar v\Big)
    &=\Re\big(\partial_t\bar v|x|^{-b}|v|^{q-2}(J_\alpha*|\cdot|^{-b}|v|^q)v\big) ,
\end{align}

which reads via the fact that ${N(q-1)+2b-\alpha-2<0}$,

\begin{align}
    \partial_t\big(\frac1{t^{N(q-1)+2b-\alpha-2}}\|\nabla v\|^2\big)
    &\geq\frac1{t^{N(q-1)+2b-\alpha-2}}\partial_t\|\nabla v\|^2\nonumber\\
    &=-\frac1{q}\partial_t\int_{\R^N}|x|^{-b}|v|^q(J_\alpha*|\cdot|^{-b}|v|^q)\,dx.\label{352}
\end{align}

Integrating in time \eqref{35}, it follows that for $0<\tau\leq t<\infty$, 

\begin{align}
   \frac1{t^{N(q-1)+2b-\alpha-2}}\|\nabla v(t)\|^2-\frac1{\tau^{N(q-1)+2b-\alpha-2}}\|\nabla v(\tau)\|^2
    &\geq-\frac1{q}\big(\int_{\R^N}|x|^{-b}|v(t)|^q(J_\alpha*|\cdot|^{-b}|v(t)|^q)\,dx\nonumber\\
    &-\int_{\R^N}|x|^{-b}|v(\tau)|^q(J_\alpha*|\cdot|^{-b}|v(\tau)|^q)\,dx\big)\,\nonumber
\end{align}

which reads for any $0<\tau\leq t<\infty$,  

\begin{align}
   \frac1{t^{N(q-1)+2b-\alpha-2}}\|\nabla v(t)\|^2+\frac1{q}\int_{\R^N}|x|^{-b}|v(t)|^q(J_\alpha*|\cdot|^{-b}|v(t)|^q)\,dx
    &\geq\frac1{\tau^{N(q-1)+2b-\alpha-2}}\|\nabla v(\tau)\|^2\nonumber\\
    &+\frac1q\int_{\R^N}|x|^{-b}|v(\tau)|^q(J_\alpha*|\cdot|^{-b}|v(\tau)|^q)\,dx\label{362}.
\end{align}

So, \eqref{36} via the mass conservation law, we get

\begin{align}
  \max\big\{\sup_{0<t\leq1}\frac1{t^{N(q-1)+2b-\alpha-2}}\|\nabla v(t)\|^2,\sup_{0<t\leq1}\|v(t)\|,\sup_{0<t\leq1}\int_{\R^N}|x|^{-b}|v(t)|^q(J_\alpha*|\cdot|^{-b}|v(t)|^q)\,dx\big\}<\infty.\label{372}
\end{align}

Take a test function $\zeta\in H^1$ and compute

\begin{align}
   \int_{\R^N} (v(t)-v(\tau))\zeta\,dx
   &=\int_\tau^t\int_{\R^N}\partial_sv(s) \zeta\,dx\,ds\nonumber\\
      &=i\int_\tau^t\int_{\R^N}\big(\Delta v-s^{N(q-1)+2b-\alpha-2}|x|^{-b}|v|^{q-2}(J_\alpha*|\cdot|^{-b}|v|^q)v \big)\zeta\,dx\,ds\nonumber\\
            &=-i\int_\tau^t\int_{\R^N}\nabla v\cdot\nabla\zeta\,dx\,ds-i\int_\tau^t\int_{\R^N}s^{N(q-1)+2b-\alpha-2}|x|^{-b}|v|^{q-2}(J_\alpha*|\cdot|^{-b}|v|^q)v\zeta\,dx\,ds\label{382}.
\end{align}

Thanks to \eqref{382} via the fact that ${N(q-1)+2b-\alpha-2>-1}$, yields

\begin{align}
    v(t)\rightharpoonup v_+,\quad\mbox{in}\quad L^2.\label{392}
\end{align}

Now, we pick $\zeta:=v(t)$ and rewrite \eqref{382} as follows

\begin{align}
  & \big|\int_{\R^N} (v(t)-v(\tau))v(t)\,dx\big|\nonumber\\
            &= \big|-i\int_\tau^t\int_{\R^N}\nabla v(s)\cdot\nabla v(t)\,dx\,ds\nonumber\\
            &-i\int_\tau^t\int_{\R^N}s^{N(q-1)+2b-\alpha-2}|x|^{-b}|v(s)|^{q-2}(J_\alpha*|\cdot|^{-b}|v(s)|^q)v(s)v(t)\,dx\,ds \big|\nonumber\\
            &\leq\int_\tau^t\|\nabla v(s)\|\|\nabla v(t)\|\,ds+\int_\tau^t\int_{\R^N}s^{N(q-1)+2b-\alpha-2}|x|^{-b}|v(s)|^{q-1}(J_\alpha*|\cdot|^{-b}|v(s)|^q)|v(t)|\,dx\,ds \label{402}.
\end{align}

Claim: there is $r_k>2$ and $0<\nu_k<\frac1q$, $k\in\{1,2\}$, such that

\begin{align}
    \big|\int_{\R^N}|x|^{-b}|v(s)|^{q-2}(J_\alpha*|\cdot|^{-b}|v(s)|^q)v(s)v(t)\,dx\big|
    &\lesssim\|v(t)\|_{r_1}^{2q-1}\|v(s)\|_{r_1}+\|v(t)\|_{r_2}^{2q-1}\|v(s)\|_{r_2}\nonumber\\
    &\lesssim\sum_{k=1}^2\big(\|\nabla v(t)\|^{2q-1}\|\nabla v(s)\|\big)^{\nu_k}.\label{clm}
\end{align}
\begin{proof}[proof of the claim]
    Using Lemma \ref{hls}, we write

\begin{align}
    (E_1)
    &:=\big|\int_{\Omega}|x|^{-b}|v(s)|^{q-2}(J_\alpha*|\cdot|^{-b}|v(s)|^q)v(s)v(t)\,dx\big|\nonumber\\
                &\lesssim \big(\||x|^{-b}\|_{L^{a_1}(\Omega)}\||x|^{-b}\|_{L^{b_1}(\Omega)}+\||x|^{-b}\|_{L^{b_2}(\Omega^c)}\||x|^{-b}\|_{L^{a_2}(\Omega)}\big)\|v(s)\|_{r_1}^{2q-1}\|v(t)\|_{r_1}\nonumber\\
                &\lesssim \|v(s)\|_{r_1}^{2q-1}\|v(t)\|_{r_1}\label{0553}.
\end{align}

Here

\begin{equation}
\left\{
\begin{array}{ll}
1+\frac{\alpha}N=\frac1{a_1}+\frac1{b_1}+\frac{2q-1}{r_1}+\frac1{r_1}=\frac1{a_2}+\frac1{b_2}+\frac{2q-1}{r_1}+\frac1{r_1} ;\\
\min\{\frac N{a_1},\frac N{a_2},\frac N{b_1}\}>{b}> \frac N{b_2};\\
r_1>2.
\label{m-13}
\end{array}
\right.
\end{equation}

So, we need

\begin{align}
    1+\frac{\alpha}N-\frac{2q}{r_1}=\frac1{a_1}+\frac1{b_1}>\frac{2b}{N}.\label{0563}
\end{align}

Then, \eqref{0563} reads

\begin{align}
     0<\frac{2q}{r_1}<1+\frac{\alpha}N-\frac{2b}{N}.\label{0573}
\end{align}

Thus, \eqref{0553} holds, since \eqref{0573} is true because 

\begin{align}
{q<1+\frac{2+\alpha-2b}{N}<1+\frac{2+\alpha-2b}{N-2}}.    
\end{align}

Also by \eqref{0573}, we have $\nu_1:=N(\frac12-\frac1{r_1})\in (0,\frac1q)$ if 

\begin{align}
    N(\frac12-\frac1{r_1})
    &<\frac1q\nonumber\\
    &\Leftrightarrow q<2(\frac q{r_1}+\frac1N)\nonumber\\
    &\Leftrightarrow q-\frac2N<\frac{2q}{r_1}.\label{ext1}
    \end{align}

    So, by \eqref{0573} and \eqref{ext1}, we need 

    \begin{align}
     q<{1+\frac{2+\alpha-2b}{N}}\label{00}.     
    \end{align}
   
    Moreover, arguing as in the above case, we get

    \begin{align}
    (E_2)
    &:=\big|\int_{\Omega^c}|x|^{-b}|v(s)|^{q-2}(J_\alpha*|\cdot|^{-b}|v(s)|^q)v(s)v(t)\,dx\big|\nonumber\\
                &\lesssim \|v(s)\|_{r_2}^{2q-1}\|v(t)\|_{r_2}\label{0553'},
\end{align}

whenever,

\begin{align}
    q>\frac{2q}{r_2}>1+\frac{\alpha}N-\frac{2b}{N}.\label{0573'}
\end{align}

Thus, \eqref{0553'} holds, since \eqref{0573'} is qossible for 

\begin{align}
{q>1+\frac{\alpha-2b}{N}}.    
\end{align}

Taking $r_2>2$, near to $2$, it follows hat $0<\nu_2<\frac1q$. The claim is proved.
    
\end{proof}
Using \eqref{372} in \eqref{402} via \eqref{clm}, we get

\begin{align}
   \big|\int_{\R^N} (v(t)-v(\tau))v(t)\,dx\big|
            &\lesssim t^{\frac{N(q-1)+2b-\alpha}{2}-1}\int_\tau^ts^{\frac{N(q-1)+2b-\alpha}{2}-1}\,ds+\sum_{k=1}^2\int_\tau^t\big(\|\nabla v(t)\|^{2q-1}\|\nabla v(s)\|\big)^{\nu_k}\,ds \nonumber\\
             &\lesssim t^{\frac{N(q-1)+2b-\alpha}{2}-1}\big(t^{\frac{N(q-1)+2b-\alpha}{2}}-\tau^{\frac{N(q-1)+2b-\alpha}{2}}\big)\nonumber\\
             &+\sum_{k=1}^2t^{(\frac{N(q-1)+2b-\alpha}{2}-1)(2q-1)\nu_k}\int_\tau^ts^{(\frac{N(q-1)+2b-\alpha}{2}-1)\nu_k}\,ds\nonumber\\
             &\lesssim t^{{N(q-1)+2b-\alpha}-1}+\sum_{k=1}^2t^{(\frac{N(q-1)+2b-\alpha}{2}-1)(2q-1)\nu_k}t^{1+(\frac{N(q-1)+2b-\alpha}{2}-1)\nu_k} .\label{0412}
\end{align}

Thus, $0<q\nu_k,t<1$ and \eqref{0412} give

\begin{align}
   \big|\int_{\R^N} (v(t)-v(\tau))v(t)\,dx\big|
             &\lesssim t^{{N(q-1)+2b-\alpha}-1}+\sum_{k=1}^2t^{({N(q-1)+2b-\alpha}-1)q\nu_k+1-q\nu_k}\nonumber\\
             &\lesssim t^{{N(q-1)+2b-\alpha}-1}+\sum_{k=1}^2t^{({N(q-1)+2b-\alpha}-1)q\nu_k}.\label{412}
\end{align}

Hence, letting $\tau\to0$ in \eqref{412} and using \eqref{392}, we write

\begin{align}
   \big|\int_{\R^N} (v(t)-v(0))v(t)\,dx\big|
             &\lesssim t^{{N(q-1)+2b-\alpha}-1}+\sum_{k=1}^2t^{({N(q-1)+2b-\alpha}-1)q\nu_k}.\label{422}
\end{align}

Now, the triangular estimate via \eqref{422} gives

\begin{align}
\|v(t)-v(0)\|^2
   &\leq\big|\int_{\R^N} (v(t)-v(0))v(t)\,dx\big|+\big|\int_{\R^N} (v(t)-v(0))v(0)\,dx\big|\nonumber\\
             &\lesssim t^{{N(q-1)+2b-\alpha}-1}+\sum_{k=1}^2t^{({N(q-1)+2b-\alpha}-1)q\nu_k}+\big|\int_{\R^N} (v(t)-v(0))v(0)\,dx\big|\nonumber\\
             &\to0,\quad\mbox{as}\quad t\to0.\label{432}
\end{align}

So, \eqref{432} via \eqref{302}, implies that for $w_+(x)=e^{it\Delta}\big((it)^{-\frac N2}e^{i\frac{|x|^2}{4t}}\bar v_0(\frac xt)\big)$,

\begin{align}
    \|w(t)-e^{-it\Delta}w_+\|
    &=\|(it)^{-\frac N2}e^{i\frac{|x|^2}{4t}}\bar v(\frac{1}{t},\frac{x}{t})-e^{-it\Delta}w_+\|\nonumber\\
        &=\|(it)^{-\frac N2}e^{i\frac{|x|^2}{4t}}\bar v(\frac{1}{t},\frac{x}{t})-(it)^{-\frac N2}e^{i\frac{|x|^2}{4t}}\bar v_0(\frac xt)\|\nonumber\\
    &=\| v(\frac{1}{t},\cdot)-v(0)\|\nonumber\\
    &\to0,\quad\mbox{as}\quad t\to\infty.\label{442}
\end{align}

The scattering is proved thanks to \eqref{442}. Theorem \ref{t2'} is established.

\hrule



\vspace{0.3cm}

{\noindent{\bf\large Declarations.}}
On behalf of all authors, the corresponding author states that there is no conflict of interest. No data-sets were generated or analyzed during the current study.

\vspace{0.3cm}

\hrule



\end{document}